\documentclass[twoside,11pt,reqno]{amsart}
\usepackage{amsmath,amssymb,amscd,mathrsfs}

\makeatletter
\@addtoreset{equation}{section}

\def\thesubsection{\thesection.\arabic{subsection}}
\def\Point{\addtocounter{subsection}{1}\vspace{2mm}
\noindent$\S$\thesubsection. \def\@currentlabel{\thesubsection}}

\newtheorem{Lemma}[equation]{Lemma}
\newtheorem{Theorem}[equation]{Theorem}
\newtheorem{Corollary}[equation]{Corollary}

\newtheorem{Example}[equation]{Example}

\def\ad{{\operatorname{ad}}}

\def\u{{\mathfrak u}}

\def\C{{\mathbb C}}

\def\Z{{\mathbb Z}}

\def\0{{\bar 0}}
\def\1{{\bar 1}}

\def\tr{{\operatorname{tr}}}
\def\str{{\operatorname{str}}}

\def\hom{{\operatorname{Hom}}}
\def\HOM{{\operatorname{\bf{Hom}}}}
\def\ext{{\operatorname{Ext}}}
\def\End{{\operatorname{End}}}

\def\soc{{\operatorname{soc}}}
\def\cosoc{{\operatorname{cosoc}}}

\def\underbar{\mathpalette\@underbar}
\def\@underbar#1#2{\settowidth{\@tempdimb}{$#1#2$}\@tempdimb=0.8\@tempdimb
                   \ooalign{$#1#2$\crcr%
                         \hfil\rule[-.5mm]{\@tempdimb}{.4pt}\hfil}}

\def\phi{{\varphi}}

\def\la{{\lambda}}

\def\De{{\Delta}}

\newdimen\hoogte    \hoogte=12pt    
\newdimen\breedte   \breedte=14pt  
\newdimen\dikte     \dikte=0.5pt 

\newenvironment{Young}{\begingroup
       \def\vr{\vrule height0.89\hoogte width\dikte depth 0.2\hoogte}
       \def\fbox##1{\vbox{\offinterlineskip
                    \hrule height\dikte
                    \hbox to \breedte{\vr\hfill##1\hfill\vr}
                    \hrule height\dikte}}
       \vbox\bgroup \offinterlineskip \tabskip=-\dikte \lineskip=-\dikte
            \halign\bgroup &\fbox{##\unskip}\unskip  \crcr }
       {\egroup\egroup\endgroup}
\def\diagram#1{\relax\ifmmode\vcenter{\,\begin{Young}#1\end{Young}\,}\else%
              $\vcenter{\,\begin{Young}#1\end{Young}\,}$\fi}

\begin{document}
\title{Tilting modules for Lie superalgebras}
\author{\sc Jonathan Brundan}
\address
{Department of Mathematics\\ University of Oregon\\
Eugene\\ OR~97403, USA}
\email{brundan@darkwing.uoregon.edu}
\thanks{Work
partially supported by the NSF (grant no. DMS-0139019).}
\maketitle

\section{Introduction}

The notion of a tilting module first emerged in Lie theory
in the 1980s, see for instance \cite{CI} where Collingwood and Irving
classified the {self-dual modules with a Verma
flag} in category $\mathcal O$ for a semisimple Lie algebra,
generalizing earlier work of Enright and Shelton \cite{ES}.
Similar looking objects were also considered by Donkin \cite{Do1}
in the representation theory of reductive algebraic groups in
positive characteristic. 
The terminology ``tilting module'' comes instead
from the representation theory of
finite dimensional algebras, via an article of Ringel \cite{Ri} 
which gives an elegant construction of tilting modules 
in the setting of quasi-hereditary algebras \cite{CPS, DR}.
Ringel's argument was subsequently
adapted to algebraic groups by Donkin \cite{Do2}
and to Lie algebras by Soergel \cite{so2}.

The goal of the present article is to extend Soergel's
framework to Lie superalgebras.
Our interest in doing this arose from the papers \cite{B1,B2}
in which we conjectured that the coefficients of certain canonical bases
should compute
multiplicities in $\De$-flags of indecomposable
tilting modules over the Lie superalgebras
$\mathfrak{gl}(m|n)$ and $\mathfrak{q}(n)$ respectively. Thus the present
article should be viewed as a companion to \cite{B1,B2}, since we provide
the general theory needed to construct the tilting modules in the first place.

We stress that the development here is very similar to Soergel's
work: most of the proofs carry over unchanged to the Lie superalgebra
setting. Like in \cite{so2}, 
we have also included in the first few sections some other well-known
generalities, most of which have their origins in the classic work of
Bernstein, Gelfand and Gelfand \cite{BGG}. The main result of the article is
best understood from Corollary~\ref{maineq2}, which roughly
speaking gives a duality
between indecomposable projective and 
indecomposable tilting modules.
The proof of this involves the construction of the ``semi-regular bimodule'',
see Lemma~\ref{icky}.

At the end of the article, we have given several examples
involving the Lie superalgebras $\mathfrak{gl}(m|n)$ and $\mathfrak{q}(n)$
to illustrate the usefulness of the theory. 
The results may also prove useful in studying the representation 
theory of the other classical Lie superalgebras and
affine Lie superalgebras.

\vspace{2mm}
\noindent
{\em Notation.}
Throughout the article, we will work over the ground field $\C$.
Suppose $V = \bigoplus_{d \in \Z} V_d = 
\bigoplus_{d \in \Z} V_{d,\0} \oplus V_{d,\1}$ is a 
{\em graded vector superspace}, i.e.
a $\Z \times \Z_2$-graded vector space.
To avoid confusion between the two different gradings, we use
the word {\em degree} to refer to the $\Z$-grading,
and {\em parity} to refer to the $\Z_2$-grading.
Write $\deg(v) \in \Z$ (resp. $\bar v \in \Z_2$) for the degree (resp.
the parity) of a homogeneous vector.
Given two graded vector superspaces $V, W$, 
$\HOM_{\C}(V, W)$ denotes the graded vector superspace with
$$
\HOM_{\C}(V, W)_{d, p} = \{f:V \rightarrow W\:|\:f(V_{d',p'}) \subseteq
W_{d+d',p+p'}\hbox{ for all }(d',p')\in\Z\times\Z_2\}
$$
for each $(d,p)\in\Z\times\Z_2$.

\section{Graded category $\mathcal O$}
For basic notions regarding Lie superalgebras, see \cite{Kac}.
Let us recall in particular that
for a Lie superalgebra $\mathfrak g = \mathfrak g_{\0}\oplus\mathfrak g_{\1}$
and $\mathfrak g$-supermodules $M, N$, a 
{homomorphism} $f:M \rightarrow N$
means a  (not necessarily even) linear map such that
$f(X m) = (-1)^{\bar f \bar X} X f(m)$ for all 
$X \in \mathfrak g, m\in M$.
This formula needs to be interpreted additively in the case that $f, X$
are not homogeneous!
We will use the notation $M \simeq N$ as opposed to the usual
$M \cong N$ to indicate that 
there is an {\em even} isomorphism between $M$ and $N$.

The category of all $\mathfrak g$-supermodules 
is not an abelian category,
but the {\em underlying even category} consisting of the same objects
and only even morphisms is abelian. This, and the existence of the
parity change functor $\Pi$, allows us to appeal to all the usual notions
of homological algebra.
Similar remarks apply to the various other categories of $\mathfrak g$-supermodules that we shall meet.

We will be concerned here instead with a
{graded} Lie superalgebra,
i.e. a Lie superalgebra $\mathfrak g$ with an additional $\Z$-grading 
$\mathfrak g = \bigoplus_{d \in \Z} \mathfrak g_d = 
\bigoplus_{d \in \Z} \mathfrak g_{d,\0}\oplus \mathfrak g_{d,\1}$
such that $[\mathfrak g_d, \mathfrak g_e] \subseteq \mathfrak g_{d+e}$
for all $d,e \in \Z$. 
A {\em graded} $\mathfrak g$-supermodule
means a $\mathfrak g$-supermodule $M$
with an additional $\Z$-grading
$M = \bigoplus_{d \in \Z} M_d = \bigoplus_{d \in \Z} M_{d,\0}\oplus M_{d,\1}$ 
such that $\mathfrak g_d M_e \subseteq M_{d+e}$ for all $d,e \in \Z$.
Homomorphisms $f:M \rightarrow N$ between graded $\mathfrak g$-supermodules
are always assumed to satisfy $f(M_d) \subseteq N_d$
for each $d \in \Z$. 

Assume from now on that we are given a 
graded Lie superalgebra 
$\mathfrak g$. Let $\mathfrak h = \mathfrak g_0,
\mathfrak b = \mathfrak g_{\geq 0} = 
\bigoplus_{d \geq 0} \mathfrak g_d$,
and $\mathfrak n = \mathfrak g_{< 0} = \bigoplus_{d < 0} \mathfrak g_d$.
We write $U(\mathfrak g), U(\mathfrak b)$ and $U(\mathfrak n)$ 
for the corresponding
universal enveloping superalgebras, all of which inherit a $\Z$-grading
from $\mathfrak g$.
We assume:
\begin{itemize}
\item[(A1)] $\dim \mathfrak g_d < \infty$ for each $d \in \Z$;
\item[(A2)] $\mathfrak h_{\0}$ is a reductive Lie algebra.
\end{itemize}
Fix in  addition a maximal toral subalgebra
$\mathfrak t$ of $\mathfrak h_{\0}$ and an abelian 
subgroup $X$ of $\mathfrak t^*$.
By an {\em admissible} representation of $\mathfrak h_{\0}$,
we mean a locally finite dimensional $\mathfrak h_{\0}$-supermodule such that
$M = \bigoplus_{\la \in X} M_\la$,
where 
$$
M_\la = \{m \in M\:|\:tm=\la(t)m\hbox{ for all }t \in \mathfrak t\}.
$$
More generally, for any graded subalgebra $\mathfrak m$ of $\mathfrak g$
containing $\mathfrak h_{\0}$, we will say that an $\mathfrak m$-supermodule
is admissible if it is admissible on restriction to $\mathfrak h_{\0}$.
We must also assume:
\begin{itemize}
\item[(A3)]
the adjoint representation $\mathfrak g$ 
is admissible.
\end{itemize}
For any graded subalgebra $\mathfrak m$ of $\mathfrak g$
containing $\mathfrak h_{\0}$, let $\mathcal C_{\mathfrak m}$ denote
the category of all admissible graded $\mathfrak m$-supermodules.
Finally let $\mathcal O$ be the category of all admissible graded
$\mathfrak g$-supermodules that are locally finite dimensional
over $\mathfrak b$.
This is a graded analogue of the category $\mathcal O$
of \cite{BGG}.

\begin{Lemma}\label{enuf}
Category $\mathcal O$ and all the categories $\mathcal C_{\mathfrak m}$
have enough injectives.
\end{Lemma}

\begin{proof}
We explain the argument for $\mathcal O$; the same argument works
for each $\mathcal C_{\mathfrak m}$.
Let $\operatorname{Fin}$ be the functor from the category of all
graded $\mathfrak g$-supermodules to $\mathcal O$ sending
an object to its largest graded submodule belonging to $\mathcal O$.
This is right adjoint to an exact functor, so sends injectives to injectives.
Moreover, the category of all graded $\mathfrak g$-supermodules has enough
injectives since it is isomorphic to the category of graded supermodules
over the universal enveloping superalgebra $U(\mathfrak g)$.
Now given any $M \in \mathcal O$, we embed
$M$ into an injective graded $\mathfrak g$-supermodule, then
apply the functor $\operatorname{Fin}$.
\end{proof}

In view of Lemma~\ref{enuf}, we can compute $\ext^i(M,N)$ in category
$\mathcal O$ or any of the categories $\mathcal C_{\mathfrak m}$
using an injective resolution of $N$.
In the sequel, we are often going to make use of the functors
$U(\mathfrak g) \otimes_{U(\mathfrak b)} ?$ and
$\HOM_{\mathfrak g_{\leq 0}}(U(\mathfrak g), ?)$.
In the latter case, for a graded $\mathfrak g_{\leq 0}$-supermodule $M$, 
$\HOM_{\mathfrak g_{\leq 0}}(U(\mathfrak g), M)$ is viewed as a 
graded $\mathfrak g$-supermodule with action $(u f)(u') = (-1)^{\bar u \bar f 
+ \bar u \bar u'} f(u'u)$, for $u, u' \in U(\mathfrak g),
f:U(\mathfrak g) \rightarrow M$.
The next lemma is a consequence of the PBW theorem.

\begin{Lemma} \label{for}
For graded
$\mathfrak b$-, $\mathfrak g_{\leq 0}$-
and $\mathfrak h$-supermodules $L, M$ and $N$,
\begin{align*}
U(\mathfrak g) \otimes_{U(\mathfrak b)} L&\simeq
U(\mathfrak g_{\leq 0}) \otimes_{U(\mathfrak h)}L,\\
U(\mathfrak g_{\leq 0}) \otimes_{U(\mathfrak h)} N&\simeq
S(\mathfrak n)\otimes N,\\\intertext{as 
graded $\mathfrak g_{\leq 0}$- resp.
$\mathfrak h$-supermodules, and}
\HOM_{\mathfrak g_{\leq 0}}(U(\mathfrak g), M)&\simeq
\HOM_{\mathfrak h}(U(\mathfrak b), M)\\
\HOM_{\mathfrak h}(U(\mathfrak b), N)&\simeq
\HOM_{\C}(S(\mathfrak g_{> 0}), N)
\end{align*}
as graded $\mathfrak b$- resp. $\mathfrak h$-supermodules.
(Here $S(\mathfrak n), S(\mathfrak g_{> 0})$ denote the symmetric
superalgebras viewed as modules via $\ad$).
\end{Lemma}

Applying the lemma and (A3),
$U(\mathfrak g) 
\otimes_{U(\mathfrak b)} ?$ (resp.
$U(\mathfrak g_{\leq 0}) \otimes_{U(\mathfrak h)} ?$)
is an exact functor from 
$\mathcal C_{\mathfrak b}$ to
$\mathcal C_{\mathfrak g}$ (resp. from
$\mathcal C_{\mathfrak h}$ to $\mathcal C_{\mathfrak g_{\leq 0}}$), 
which is obviously left adjoint to the natural
restriction functor. Similarly, 
$\HOM_{\mathfrak g_{\leq 0}}(U(\mathfrak g), ?)$ is an exact functor
from $\mathcal C_{\mathfrak g_{\leq 0}}$ to $\mathcal C_{\mathfrak g}$
that is right adjoint to restriction.

\begin{Lemma}\label{gfr}
For $i \geq 0$, $L \in \mathcal C_{\mathfrak g}$,
$M \in \mathcal C_{\mathfrak h}$ and
$N \in \mathcal C_{\mathfrak g_{\leq 0}}$, we have that
\begin{align*}
\ext^i_{\mathcal C_{\mathfrak g}}(L, \HOM_{\mathfrak g_{\leq 0}}(U(\mathfrak g), N))
&\simeq \ext^i_{\mathcal C_{\mathfrak g_{\leq 0}}}(L, N),\\
\ext^i_{\mathcal C_{\mathfrak g_{\leq 0}}}(U(\mathfrak g_{\leq 0}) \otimes_{U(\mathfrak h)} M, N) &\simeq
\ext^i_{\mathcal C_{\mathfrak h}}(M, N).
\end{align*}
\end{Lemma}

\begin{proof}
Argue by induction on $i$ using the long exact sequence.
\end{proof}

\section{Standard and costandard modules}
Let $\Lambda$ be a complete set of pairwise
non-isomorphic irreducible admissible graded $\mathfrak h$-supermodules.
Each $E \in \Lambda$ is necessarily concentrated in a single degree,
denoted $|E| \in \Z$. 
Moreover, by the superalgebra analogue of Schur's lemma,
the number
\begin{equation}
d_E := \dim \End_{\mathcal C_{\mathfrak h}}(E)
\end{equation}
is either $1$ or $2$.

\begin{Lemma}\label{tobegin}
Every $E \in \Lambda$
has a finite dimensional projective cover $\widehat E$ in
$\mathcal C_{\mathfrak h}$, with $\cosoc_{\mathfrak h} \widehat E \simeq E$.
Moreover, given objects $M, P \in \mathcal C_{\mathfrak h}$ with
$P$ projective, $M \otimes P$ is also projective.
\end{Lemma}

\begin{proof}
For a graded $\mathfrak h_{\0}$-supermodule $M$, we observe that
\begin{equation*}
U(\mathfrak h) \otimes_{U(\mathfrak h_{\0})} M
\simeq S(\mathfrak h_{\1}) \otimes M
\end{equation*}
as graded $\mathfrak h_{\0}$-supermodules.
Combining this with (A3) shows that
the functor
$U(\mathfrak h) \otimes_{U(\mathfrak h_{\0})} ?$ maps
$\mathcal C_{\mathfrak h_{\0}}$ to $\mathcal C_{\mathfrak h}$. 
Since it is left adjoint to an exact functor, it maps
projectives to projectives. By (A2) and Weyl's theorem on
complete reducibility, every object in $\mathcal C(\mathfrak h_{\0})$ 
is projective.

Now take $E \in \Lambda$.
Let $\widehat E$ be any indecomposable summand of 
$U(\mathfrak h) \otimes_{U(\mathfrak h_{\0})} E$
that maps surjectively onto $E$ under the natural multiplication map.
By the preceeding paragraph, $\widehat E$ is a finite dimensional 
indecomposable projective object in $\mathcal C_{\mathfrak h}$ 
mapping surjectively
onto $E$. Now the usual arguments via Fitting's lemma show that
$\widehat E$ is actually a projective cover of $E$ 
in the category $\mathcal C_{\mathfrak h}$ and that
$\cosoc_{\mathfrak h} \widehat{E} \simeq E$.

Finally let $P \in \mathcal C_{\mathfrak h}$ 
be an arbitrary projective object. 
Then, we can find
$Q \in \mathcal C_{\mathfrak h}$ and $R \in \mathcal C_{\mathfrak h_\0}$
such that $P \oplus Q \cong U(\mathfrak h) \otimes_{U(\mathfrak h_{\0})}
R$. By the tensor identity, $(P \oplus Q) \otimes M
\cong U(\mathfrak h) \otimes_{U(\mathfrak h_{\0})} (R \otimes M)$.
The latter is projective and $P \otimes M$ is isomorphic to 
a summand of it, so 
$P \otimes M$ is projective too.
\end{proof}

Define the {\em standard} and {\em costandard $\mathfrak g$-supermodules} corresponding to $E\in \Lambda$:
\begin{equation}
\Delta(E) := U(\mathfrak g) \otimes_{U(\mathfrak b)} \widehat{E},
\qquad
\nabla(E) := \HOM_{\mathfrak g_{\leq 0}}(U(\mathfrak g), E).
\end{equation}
By Lemma~\ref{for}, both $\Delta(E)$ and $\nabla(E)$ 
are admissible, and 
clearly they are locally finite dimensional over $\mathfrak b$, hence
they belong to $\mathcal O$.
Indeed, 
letting $\mathcal O_{\leq d}$ denote the full subcategory of $\mathcal O$
consisting of all objects that are zero in degrees $> d$,
both $\Delta(E)$ and $\nabla(E)$ belong to $\mathcal O_{\leq |E|}$,
with $\Delta(E)_{|E|} \simeq \widehat{E}$, $\nabla(E)_{|E|} \simeq E$.
We define
\begin{equation}
L(E) := \cosoc_{\mathfrak g} \Delta(E)
\end{equation}
for each $E \in \Lambda$.
The following well-known lemma shows in particular that these are
irreducible.

\begin{Lemma}\label{class} 
The $\{L(E)\}_{E \in \Lambda}$ form a complete
set of pairwise non-isomorphic irreducibles in $\mathcal O$.
Moreover, $L(E) \simeq \operatorname{soc}_{\mathfrak g} \nabla(E)$.
\end{Lemma}

\begin{proof}
Over $\mathfrak g_{\leq 0}$, 
$\Delta(E) \simeq U(\mathfrak g_{\leq 0}) \otimes_{U(\mathfrak h)} \widehat{E}$,
hence $\cosoc_{\mathfrak g_{\leq 0}} \Delta(E)
\simeq \cosoc_{\mathfrak h} \widehat{E} \simeq E$.
This immediately implies that $L(E)$ is irreducible in $\mathcal O$
and $\cosoc_{\mathfrak g_{\leq 0}} L(E) \simeq E$.
Hence the $\{L(E)\}_{E \in \Lambda}$ are pairwise non-isomorphic
irreducibles. Now take any irreducible $M \in \mathcal O$.
There exists a non-zero $\mathfrak b$-homomorphism $E \rightarrow M$
for some $E \in \Lambda$. This induces by Frobenius reciprocity
a non-zero $\mathfrak g$-homomorphism $\Delta(E) \rightarrow M$, 
hence $M \cong L(E)$. The same argument shows that $\soc_{\mathfrak b}
L(E) \simeq E$. Finally, over $\mathfrak b$,
$\nabla(E) \simeq \HOM_{\mathfrak h}(U(\mathfrak b), E)$,
so $\soc_{\mathfrak b} \nabla(E) \simeq E$. Hence $\soc_{\mathfrak g}
\nabla(E) \simeq L(E)$ too.
\end{proof}

\begin{Lemma}\label{sl1} 
Let $E, F \in \Lambda$.
\begin{itemize}
\item[(i)] $\Delta(E)$ is the projective cover of $L(E)$ in
$\mathcal O_{\leq |E|}$.
\item[(ii)]
$\dim \hom_{\mathcal O}(\Delta(E), \nabla(F)) = 0$
if $E \neq F$, $d_E$ if $E = F$.
\item[(iii)] $\ext^1_{\mathcal O}(\Delta(E), \nabla(F)) = 0$.
\end{itemize}
\end{Lemma}

\begin{proof}
For (i), take $M \in \mathcal O_{\leq |E|}$. We have the following
sequence of isomorphisms natural in $M$:
$$
\hom_{\mathcal O_{\leq |E|}}(\Delta(E), M)
\simeq \hom_{\mathcal C_{\mathfrak g}}(\Delta(E), M)
\simeq \hom_{\mathcal C_{\mathfrak b}}(\widehat{E}, M)
\simeq \hom_{\mathcal C_{\mathfrak h}}(\widehat{E}, M).
$$
Since $\widehat{E}$ is projective in $\mathcal C_{\mathfrak h}$,
this shows that $\Delta(E)$ is projective in $\mathcal O_{\leq |E|}$.
The same argument with $M = \Delta(E)$ shows that
$\dim \End_{\mathcal O_{\leq|E|}}(\Delta(E))$ is finite dimensional,
so we get that $\Delta(E)$ is actually the projective cover of $L(E)$ in 
$\mathcal O_{\leq |E|}$ from Fitting's lemma.
For
(ii), (iii), Lemma~\ref{gfr} implies 
for every $i \geq 0$ 
that
$$
\ext^i_{\mathcal C_{\mathfrak g}}(\Delta(E), \nabla(F))
\simeq
\ext^i_{\mathcal C_{\mathfrak h}}(\widehat{E}, F).
$$
Since $\widehat{E}$ is projective with $\cosoc_{\mathfrak h} \widehat{E} \simeq E$,
the right hand side is zero if $i > 0$ or if $E \neq F$,
and is of dimension $d_E$ otherwise.
Now we are done since $\mathcal O$ is a full subcategory of
$\mathcal C_{\mathfrak g}$.
\end{proof}

\section{Projective modules and blocks}
Let $M \in \mathcal O$.
A {\em $\Delta$-flag} of $M$ means a filtration
$$
0 = M_0 \subseteq M_1 \subseteq M_2 \dots
$$
such that $M = \bigcup_{i\geq 0} M_i$ and each factor
$M_i / M_{i-1}$ is either zero or $\cong \Delta(E_i)$
for $E_i \in \Lambda$.
If the filtration stabilizes after finitely many terms
we will call it a {\em finite $\Delta$-flag}.
Arguing as in \cite[Lemma 5.10]{so2}, one shows:

\begin{Lemma}\label{indstep}
Suppose we have that $\ext^1_{\mathcal O}(\Delta(F), N) = 0$ for all
$F \in \Lambda$.
Then,
$\ext^1_{\mathcal O}(M, N) = 0$ for every $M \in \mathcal O$
admitting a $\Delta$-flag.
\end{Lemma}

Applying the lemma to $N = \nabla(E)$, one easily deduces
that the multiplicity
of $\Delta(E)$ as a subquotient of a $\Delta$-flag
of $M$ is equal to 
$\dim \hom_{\mathcal O}(M, \nabla(E)) / d_E$,
for every $M \in \mathcal O$ admitting a $\Delta$-flag.
In particular, this multiplicity does not depend on the choice of the
$\Delta$-flag. We will denote it by $(M:\Delta(E))$.

\begin{Lemma}\label{dep}
A graded $\mathfrak g$-supermodule $M$ admits a finite $\Delta$-flag
if and only if 
$M$ is a graded free $U(\mathfrak n)$-supermodule
of finite rank and its restriction to $\mathfrak h$ is a projective
object in $\mathcal C_{\mathfrak h}$.
\end{Lemma}

\begin{proof}
($\Rightarrow$) It suffices to prove this for $M = \Delta(E)$.
Obviously this is a graded free $U(\mathfrak n)$-supermodule
of rank $\dim \widehat{E}$. Moreover, over $\mathfrak h$,
we have by Lemma~\ref{for} that
$M \simeq S(\mathfrak n) \otimes \widehat{E}$.
This is projective in $\mathcal C_{\mathfrak h}$ by Lemma~\ref{tobegin}.

($\Leftarrow$) 
We may assume that $M = \bigoplus_{i=1}^n U(\mathfrak n) \otimes V_i$
is a decomposition of $M$ as a graded free $U(\mathfrak n)$-supermodule,
where $V_i$ is a finite dimensional vector superspace 
concentrated in degree $d_i$ with trivial action of $\mathfrak n$,
and $d_1 > \dots > d_n$.
Note then that $1 \otimes V_1$ must be invariant under the action
of $\mathfrak b$, and $\mathfrak g_{> 0}$ acts trivially.
Hence by the projectivity assumption
it decomposes as a direct sum of finitely many $\widehat{E}$'s as a 
$\mathfrak b$-supermodule. Each $U(\mathfrak n) \otimes \widehat{E}$
in this decomposition is isomorphic as a graded $\mathfrak g$-supermodule
to $\Delta(E)$, and the quotient of $M$ by $U(\mathfrak n) \otimes V_1$
is graded free of strictly smaller rank and is still
projective over $\mathfrak h$, so we are done by induction.
\end{proof}

\begin{Corollary}\label{fp} If $M$ admits a finite $\Delta$-flag, so
does any summand of $M$.
\end{Corollary}

\begin{proof}
Any summand of
a graded free $U(\mathfrak n)$-supermodule of finite rank
is again graded free of finite rank, see
\cite[Remark 2.4(2)]{so2}.
\end{proof}

We now come to the 
basic result on projective objects in category $\mathcal O$.

\begin{Theorem}\label{pct}
Every simple object $L(E) \in \mathcal O_{\leq n}$
admits a projective cover $P_{\leq n}(E)$ in $\mathcal O_{\leq n}$
with $\cosoc_{\mathfrak g} P_{\leq n}(E) \simeq L(E)$.
Moreover,
\begin{itemize}
\item[(i)] $P_{\leq n}(E)$ admits 
a finite $\Delta$-flag with $\De(E)$ at the top;
\item[(ii)] for $m > n$, the kernel of any surjection
$P_{\leq m}(E) \twoheadrightarrow P_{\leq n}(E)$ admits a finite $\Delta$-flag
with subquotients of the form $\Delta(F)$ for $m \geq |F| > n$;
\item[(iii)] $L(E)$ admits a projective cover $P(E)$ in $\mathcal O$
if and only if there exists $n \gg 0$ 
with
$P_{\leq n}(E) = P_{\leq n+1}(E) = \dots$, 
in which case $P(E) = P_{\leq n}(E)$.
\end{itemize}
\end{Theorem}

\begin{proof}
The proof is essentially the same as \cite[Theorem 3.2]{so2}, so we
just sketch the construction of $P_{\leq n}(E)$ and refer the reader to
{\em loc. cit.} for everything else.
For a graded $\mathfrak b$-supermodule $M$, let $\tau_{\leq n} M$
denote the quotient of $M$
by the submodule $\bigoplus_{d > n} M_d$ 
of all homomogeneous parts of degree $> n$.
For $E \in \Lambda$,
$$
Q := U(\mathfrak g) \otimes_{U(\mathfrak b)} \tau_{\leq n}(U(\mathfrak b) \otimes_{U(\mathfrak h)} \widehat{E})
$$
is projective in $\mathcal O_{\leq n}$ as in the proof
of \cite[Theorem 3.2(1)]{so2}, it is graded free over 
$U(\mathfrak n)$ of finite rank, and it is 
projective viewed as an object of $\mathcal C_{\mathfrak h}$
by Lemma~\ref{tobegin}.
So Lemma~\ref{dep} shows that 
$Q$ has a finite $\Delta$-flag.
Now $Q$ clearly maps surjectively onto $L(E)$.
Let $P_{\leq n}(E)$ be an indecomposable summand of $Q$ that also
maps surjectively onto $L(E)$. This has a finite $\Delta$-flag too
by Corollary~\ref{fp}, and it is a projective cover of $L(E)$
in $\mathcal O_{\leq n}$ by a Fitting's lemma argument, see
\cite[Lemma 3.3]{so2}.
\end{proof}

\vspace{1mm}

For $M \in \mathcal O$, we 
write $[M:L(E)]$ for the composition multiplicity of $L(E)$
in $M$, i.e. the supremum
of
$\#\{i\:|\:M_i / M_{i-1} \cong L(E)\}$
over all finite filtrations
$M = (M_i)_i$ of $M$.
This multiplicity is additive on short exact sequences.
Now we get ``BGG reciprocity'':

\begin{Corollary}\label{bggrec}
$(P_{\leq n}(E) : \Delta(F))
= [\nabla(F):L(E)]$
for all $E, F \in \Lambda$ and $n \geq |E|, |F|$.
\end{Corollary}

\begin{proof}
In $\mathcal O_{\leq n}$, we have that
$[\nabla(F):L(E)] = \dim \hom_{\mathcal O}(P_{\leq n}(E), \nabla(F)) / d_E$.
This equals $(P_{\leq n}(E):\Delta(F))$ by the definition of the latter
multiplicity.
\end{proof}

Suppose finally in this section that $\sim$ is an equivalence relation
on $\Lambda$ with the property that
\begin{equation*}
[\Delta(F):L(E)] \neq 0 \hbox{ or }
[\nabla(F):L(E)] \neq 0 \Rightarrow F \sim E
\end{equation*}
for each $E, F \in \Lambda$.
For an equivalence class $\theta \in \Lambda / \sim$, let
$\mathcal O_\theta$ be the full subcategory of $\mathcal O$
consisting of the objects $M \in \mathcal O$ all of whose
irreducible subquotients are of the form
$L(E)$ for $E \in \theta$.
We refer to $\mathcal O_\theta$ as a {\em block}
of $\mathcal O$, in view of the following theorem which
is proved exactly as in
\cite[Theorem 4.2]{so2}.

\begin{Theorem}\label{blockdec}
The functor $$
\prod_{\theta \in \Lambda / \sim}
\mathcal O_\theta \rightarrow \mathcal O,\quad
(M_\theta)_{\theta} \mapsto \bigoplus_{\theta \in \Lambda / \sim} M_\theta
$$
is an equivalence of categories.
\end{Theorem}

\section{Tilting modules and Arkhipov-Soergel duality}
Next, we discuss the classification of 
{tilting modules} in $\mathcal O$.
The first main result is the analogue of \cite[Theorem 5.2]{so2}.

\begin{Theorem}\label{tiltthm} For any $E \in \Lambda$, 
there exists a unique up to isomorphism indecomposable object 
$T(E) \in \mathcal O$ such that
\begin{itemize}
\item[(i)] $\ext^1_{\mathcal O}(\Delta(F), T(E)) = 0$ for all 
$F \in \Lambda$;
\item[(ii)] $T(E)$ admits a $\Delta$-flag starting with $\Delta(E)$ at the 
bottom.
\end{itemize}
\end{Theorem}

We call $T(E)$ the {\em indecomposable tilting module} corresponding to 
$E \in \Lambda$.
The proof given by Soergel is a variation on an argument of Ringel
\cite{Ri}, and carries over to the present setting virtually unchanged.
The main step is to show that for any $E \in \Lambda$ with $|E| \geq n$, there
exists a unique up to isomorphism
indecomposable object 
$T_{\geq n}(E)$ in $\mathcal O$
such that 
\begin{itemize}
\item[(i)$'$] $\ext^1_{\mathcal O}(\Delta(F), T_{\geq n}(E)) = 0$
for all $F \in \Lambda$ with 
$|F| \geq n$;
\item[(ii)$'$]
$T_{\geq n}(E)$ admits a finite $\Delta$-flag starting
with $\Delta(E)$ at the bottom and with all other subquotients
of the form $\Delta(F)$ for $F$'s 
with $|E| > |F| \geq n$.
\end{itemize}
Moreover, given $|E| \geq m \geq n$, there exists an inclusion
$T_{\geq m}(E) \hookrightarrow T_{\geq n}(E)$, and the cokernel 
of any such inclusion admits a 
finite $\Delta$-flag with subquotients $\Delta(F)$ for $m > |F| \geq n$.
Given these results, a candidate for the 
desired module $T(E)$ can then be constructed 
as a direct limit of the $T_{\geq n}(E)$'s as $n \rightarrow -\infty$.
Uniqueness then needs to be established separately.

To proceed, we need to make two additional assumptions
(see \cite[Remark 1.2]{so2} for remarks on the first one):
\begin{itemize}
\item[(A4)] $\mathfrak g$ is generated as a Lie superalgebra by
$\mathfrak g_0, \mathfrak g_1$ and $\mathfrak g_{-1}$;
\item[(A5)] for $E \in \Lambda$, $(\widehat E)^* \cong \widehat{E^\#}$
for some $E^\# \in \Lambda$.
\end{itemize}
Under the assumption (A4), 
an {\em admissible semi-infinite character} $\gamma$ for $\mathfrak g$
is defined to be a Lie superalgebra homomorphism
$\gamma:\mathfrak h \rightarrow \C$ such that
$\gamma|_{\mathfrak t}\in X$ and
\begin{equation}
\gamma([X, Y]) = \str_{\mathfrak h} (\ad X \circ \ad Y)
\end{equation}
for all $X \in \mathfrak g_1, Y \in \mathfrak g_{-1}$.
(We recall the {\em supertrace} of an endomorphism
$f=f_{\0}+f_{\1}:V \rightarrow V$ of a vector superspace is defined
by $\str_V f := \tr_{V_{\0}} f_{\0} - \tr_{V_{\1}} f_{\0}$.)

In the next lemma, we 
write $U(\mathfrak n)^{\circledast}$ for the
graded dual $\HOM_{\C}(U(\mathfrak n), \C)$ (where $\C = \C_{0,\0}$) 
viewed as a $U(\mathfrak n),U(\mathfrak n)$-bimodule with
left and right actions defined by $(nf)(n') = (-1)^{\bar n \bar f + 
\bar n \bar n'} f(n'n)$ and $(fn)(n') = f(nn')$ respectively,
for $n, n' \in U(\mathfrak n), f \in U(\mathfrak n)^\circledast$.

\begin{Lemma}\label{icky}
Let $\gamma:\mathfrak h \rightarrow \C$ be an admissible 
semi-infinite character
for $\mathfrak g$.
Then there exists a graded $U(\mathfrak g), U(\mathfrak g)$-bimodule 
$S_\gamma$ and an even monomorphism
$\iota:U(\mathfrak n)^\circledast \hookrightarrow S_\gamma$ 
of graded $U(\mathfrak n),U(\mathfrak n)$-bimodules such that
\begin{itemize}
\item[(i)] the map $U(\mathfrak g) 
\otimes_{U(\mathfrak n)} U(\mathfrak n)^\circledast \rightarrow S_\gamma,
u \otimes f \mapsto u \iota(f)$ is a bijection;
\item[(ii)] the map $U(\mathfrak n)^\circledast 
\otimes_{U(\mathfrak n)} U(\mathfrak g) \rightarrow S_\gamma,
f \otimes u \mapsto \iota(f) u$ is a bijection;
\item[(iii)] 
$[H, \iota(f)] 
= \iota(f) \gamma(H) - (-1)^{\bar H \bar f} \iota(f \circ \ad H)$
for all $H \in \mathfrak h$ and $f \in U(\mathfrak n)^\circledast$.
\end{itemize}
\end{Lemma}

\begin{proof}
This is proved in almost exactly the same way as \cite[Theorem 1.3]{so2}.
However, the signs are rather delicate in the super case. So we 
describe explicitly the construction of $S_\gamma$, referring to
the proof of \cite[Theorem 1.3]{so2} for a fuller account of the other steps
that need to be made.
As a graded vector superspace, we have that
$$
S_\gamma = U(\mathfrak n)^\circledast \otimes_{\C} U(\mathfrak b),
$$
and the map
$\iota:U(\mathfrak n)^\circledast \rightarrow S_\gamma$ is defined by
$\iota(f) = f\otimes 1$.
Note $S_\gamma$ is a $U(\mathfrak n), U(\mathfrak b)$-bimodule in the usual
way.
We now extend this structure
to make $S_\gamma$ into $U(\mathfrak g), U(\mathfrak g)$-bimodule.
First, there is a natural isomorphism of 
$U(\mathfrak n), U(\mathfrak b)$-bimodules
$$
S_\gamma = U(\mathfrak n)^\circledast \otimes_{\C} U(\mathfrak b)
\stackrel{\sim}{\longrightarrow}
U(\mathfrak n)^\circledast \otimes_{U(\mathfrak n)} 
U(\mathfrak g)
$$
mapping $u \otimes v$ to $u \otimes v$; we get the right action of $U(\mathfrak g)$
on $S_\gamma$ via this isomorphism.
To obtain the left action, we
use the natural isomorphisms
$$
S_\gamma = U(\mathfrak n)^\circledast \otimes_{\C} U(\mathfrak b)
\stackrel{\sim}{\longrightarrow}
\HOM_{\C}(U(\mathfrak n), U(\mathfrak b))
\stackrel{\sim}{\longleftarrow}
\HOM_{U(\mathfrak b)}(U(\mathfrak g), \C_\gamma \otimes_\C U(\mathfrak b)).
$$
For the right hand space,
the action of $U(\mathfrak b)$ 
is the natural left action on $U(\mathfrak g)$, and the tensor product
of the action on $\C_\gamma = \C_{0,\0}$ affording the character $\gamma$
and the natural left action on $U(\mathfrak b)$.
The first isomorphism maps $f \otimes b$ to the function
$\widehat{f \otimes b}:n \mapsto (-1)^{\bar b \bar n}f(n) b$. 
The second isomorphism
is given by restriction of functions from 
$U(\mathfrak g)$ to $U(\mathfrak n)$, identifying
$\C_\gamma \otimes_{\C} U(\mathfrak b)$ with
$U(\mathfrak b)$ via $1 \otimes u \mapsto u$.
Now, $U(\mathfrak g)$ acts naturally on the left on the right hand space,
by $(uf)(u') = (-1)^{\bar u \bar f + \bar u \bar u'} f(u'u)$,
for $u,u' \in U(\mathfrak g)$ and $f:U(\mathfrak g) \rightarrow \C_\gamma \otimes_{\C} U(\mathfrak b)$.
Transferring this to $S_\gamma$ via the isomorphisms gives
the left $U(\mathfrak g)$-module structure on $S_\gamma$.
Now we have to check that the left and right actions
of $U(\mathfrak g)$ on $S_\gamma$ just defined
commute with one another, so that $S_\gamma$ is a $U(\mathfrak g), U(\mathfrak g)$-bimodule. This is done by brutal calculation relying on the
assumption that $\gamma$ is a semi-infinite character, see
the proof of \cite[Theorem 1.3]{so2} for the detailed argument which
generalizes routinely to our setting. Once that is done, (i)--(iii) are
relatively easy to check to complete the proof.
\end{proof}

For the remainder of the section, we fix an admissible 
semi-infinite character
$\gamma$ for $\mathfrak g$
and let $S_\gamma$ be the {\em semi-regular
bimodule}
constructed in Lemma~\ref{icky}.
Let $\mathcal M$ resp. $\mathcal K$
be the category of all admissible graded $\mathfrak g$-supermodules
that are free resp. cofree of finite rank as graded 
$U(\mathfrak n)$-supermodules, i.e. isomorphic to direct sums of
maybe graded shifted copies of $U(\mathfrak n)$ resp. $U(\mathfrak n)^{\circledast}$.
The following theorem is the super analogue of
\cite[Theorem 2.1]{so2}, which Soergel attributes originally
to Arkhipov \cite{Ar}.

\begin{Theorem}\label{at}
The functors $\mathcal M \rightarrow \mathcal K,
M \mapsto S_\gamma \otimes_{U(\mathfrak g)} M$
and $\mathcal K \rightarrow \mathcal M,
M \mapsto \HOM_{U(\mathfrak g)}(S_\gamma, M)$ are mutually
inverse equivalences between the categories $\mathcal M$ and $\mathcal K$,
such that short exact sequences correspond to short exact sequences.
\end{Theorem}

\begin{proof}
Take $M \in \mathcal M$.
Recalling Lemma~\ref{icky},
the map $f \otimes m \mapsto \iota(f) \otimes m$
is a $U(\mathfrak n)$-isomorphism
$U(\mathfrak n)^\circledast \otimes_{U(\mathfrak n)} M
\rightarrow S_\gamma \otimes_{U(\mathfrak g)} M$.
Hence $S_\gamma \otimes_{U(\mathfrak g)} M$ is graded cofree
of finite rank, so in particular it is finite dimensional in each degree.
Moreover, for $f \in U(\mathfrak n)^\circledast, m \in M$ and
$H \in \mathfrak h$, we have by Lemma~\ref{icky}(iii) that
\begin{equation}\label{act}
H (\iota(f) \otimes m) = (-1)^{\bar H \bar f}
\iota(f) \otimes (H + \gamma(H)) m
- (-1)^{\bar H \bar f} \iota(f \circ \ad H) \otimes m.
\end{equation}
It follows from this and (A3) that $S_\gamma \otimes_{U(\mathfrak g)}
M$ is admissible. Hence $S_\gamma \otimes_{U(\mathfrak g)} ?$
is a well-defined functor from $\mathcal M$ to $\mathcal K$.
For the other direction, 
we  note that $\HOM_{U(\mathfrak n)}(U(\mathfrak n)^\circledast,
U(\mathfrak n)^{\circledast}) \simeq U(\mathfrak n)$ as a 
$U(\mathfrak n), U(\mathfrak n)$-bimodule; an isomorphism maps
$u \in U(\mathfrak n)$ to 
$\widehat u \in \HOM_{U(\mathfrak n)}(U(\mathfrak n)^\circledast,
U(\mathfrak n)^{\circledast})$ where
$(\widehat u f)(n) = (-1)^{\bar u \bar f} f(un)$ for each
$f \in U(\mathfrak n)^{\circledast}, n \in U(\mathfrak n)$.
So for $N \in \mathcal K$, we deduce that
$\HOM_{U(\mathfrak g)}(S_\gamma, N) \simeq
\HOM_{U(\mathfrak n)}(U(\mathfrak n)^{\circledast}, N)$ is 
graded free of finite rank over $U(\mathfrak n)$.
Moreover, given $\theta \in \HOM_{U(\mathfrak g)}(S_\gamma, N)$, 
\begin{equation}\label{tip}
(H \theta)(\iota(f)) = 
(H-\gamma(H))\theta(\iota(f))
+
(-1)^{\bar H \bar \theta+\bar H \bar f}
\theta(\iota(f \circ \ad H))
\end{equation}
for each $H \in \mathfrak h$ and 
$f \in U(\mathfrak n)^{\circledast}$.
Using this and (A3) one can check that
$\HOM_{U(\mathfrak g)}(S_\gamma, N)$ is admissible.
Hence, $\HOM_{U(\mathfrak g)}(S_\gamma, ?)$
is a well-defined functor from $\mathcal K$ to $\mathcal M$.
The remainder of the proof 
is exactly as in the proof of \cite[Theorem 2.1]{so2}.
\end{proof}

Finally, let $\mathcal O^\Delta$ be the full subcategory of $\mathcal O$ 
consisting of all objects admitting a finite $\Delta$-flag.
We recall from Corollary~\ref{fp} 
that $\mathcal O^\Delta$ is closed under
taking direct summands.
For a graded $\mathfrak g$-supermodule $M$, we let
$M^\star$ denote its graded dual, namely, the space
$\HOM_{\C}(M, \C)$, where $\C = \C_{0,\0}$, with action defined by
$(Xf)(m) = -(-1)^{\bar X \bar f}f(X m)$ for each $X \in \mathfrak g,
m \in M$ and $f:M \rightarrow \C$.
Recalling the assumption (A5), the theorem has the following 
corollary:

\begin{Corollary}\label{maineq2}
The functor $M \mapsto (S_\gamma \otimes_{U(\mathfrak g)} M)^\star$
defines a contravariant 
equivalence of categories $\mathcal O^\Delta \rightarrow
\mathcal O^\Delta$
under which short exact sequences correspond to short exact sequences, 
$\Delta(\C_{-\gamma} \otimes E^\#)$ maps 
to $\Delta(E)$ 
and $P_{\leq -n}(\C_{-\gamma} \otimes E^\#)$ maps
to $T_{\geq n}(E)$,
for every $E \in \Lambda$ and $n \leq |E|$.
\end{Corollary}

\begin{proof}
It is easy to see using (\ref{act}) and (A5) that the degree $-|E|$ 
piece of
$(S_\gamma \otimes_{U(\mathfrak g)} \Delta(E))^\star
\simeq (U(\mathfrak n)^\circledast \otimes \widehat{E})^\star$
is isomorphic to $\C_{-\gamma} \otimes 
\widehat{E^\#}$ as an $\mathfrak h$-supermodule.
Moreover, this generates
$(S_\gamma \otimes_{U(\mathfrak g)} \Delta(E))^\star$ freely
as a $U(\mathfrak n)$-supermodule, hence
$(S_\gamma \otimes_{U(\mathfrak g)} \Delta(E))^\star
\cong \Delta(\C_{-\gamma} \otimes E^\#)$.
It follows from this and Theorem~\ref{at}
that the functor $(S_\gamma \otimes_{U(\mathfrak g)} ?)^\star$
maps $\mathcal O^\Delta$ to $\mathcal O^\Delta$ and sends
short exact sequences to short exact sequences.
Similarly, one shows using (\ref{tip}) that
$\HOM_{U(\mathfrak g)}(S_\gamma, \Delta(E)^\star) \cong
\Delta(\C_{-\gamma} \otimes E^\#)$.
Hence the functor $\HOM_{U(\mathfrak g)}(S_\gamma, ?^\star)$
maps $\mathcal O^\Delta$ to $\mathcal O^\Delta$.
Now it is immediate from Theorem~\ref{at} that our two functors
are mutually inverse equivalences.
It just remains to show that $(S_\gamma \otimes_{U(\mathfrak g)}
P_{\leq -n}(\C_{-\gamma} \otimes E^\#))^\star \cong T_{\geq n}(E)$,
for $n \leq |E|$, for which one uses the characterization
of $T_{\geq n}(E)$ given in (i)$'$, (ii)$'$ above.
\end{proof}

\begin{Corollary}\label{dde}
For $E, F \in \Lambda$, we have that
$$
(T(E):\Delta(F)) = [\nabla(\C_{-\gamma} \otimes F^\#) : 
L(\C_{-\gamma} \otimes E^\#)].
$$
\end{Corollary}

\begin{proof}
We have for $n \leq |E|, |F|$ that 
\begin{align*}
(T(E):\Delta(F)) &= (T_{\geq n}(E): \Delta(F))
=
(P_{\leq -n}(\C_{-\gamma} \otimes E^\#):\Delta(\C_{-\gamma}\otimes F^\#))\\
&= [\nabla(\C_{-\gamma} \otimes F^\#):L(\C_{-\gamma} \otimes E^\#)],
\end{align*}
using Corollary~\ref{maineq2} and Lemma~\ref{bggrec}.
\end{proof}

\section{Some variations}\label{var}
We now mention some variations to the general framework considered so far.
First of all, we recall from
\cite[$\S$6]{so2} how to deduce results about ungraded
$\mathfrak g$-supermodules from the graded theory above.
To do this, one needs to require in addition that
\begin{itemize}
\item[(A6)] there is an element $D \in \mathfrak h_{\0}$ such that
$[D, X] = \deg(X) X$ for all homogeneous $X \in \mathfrak g$.
\end{itemize}
Let $\overline{\mathcal O}$ be the category of all
admissible (but no longer graded!) 
$\mathfrak g$-supermodules that are 
locally finite dimensional over $\mathfrak b$.
Since $D$ necessarily belongs to $\mathfrak t$,
every $M \in {\mathcal O}$ resp. $M\in \overline{\mathcal O}$
decomposes into eigenspaces
$M = \bigoplus_{a \in \C} M^{(a)}$
with respect to the action of $D$.
For $a \in \C$, let $\mathcal O_a$ denote the full subcategory
of $\mathcal O$ consisting of all $M\in\mathcal O$ such that $M^{(a+i)} = M_i$ for all $i \in \Z$.
For $\bar a \in \C / \Z$, let $\overline{\mathcal O}_{\bar a}$
denote the full subcategory of $\overline{\mathcal O}$
consisting of all $M\in\overline{\mathcal O}$
such that $M^{(b)} = 0$ for all $b \notin \bar a$.
Then,
\begin{equation*}\label{bldec}
\mathcal O = \prod_{a \in \C} \mathcal O_a,
\qquad
\overline{\mathcal O} = \prod_{\bar a \in \C / \Z}
\overline{\mathcal O}_{\bar a}.
\end{equation*}
Forgetting the grading gives an isomorphism of categories
$\mathcal O_a \rightarrow \overline{\mathcal O}_{\bar a}$,
the inverse functor being defined on $M \in \overline{\mathcal O}_{\bar a}$
by introducing a $\Z$-grading
according to the rule $M_i = M^{(a+i)}$.
In this way, we can transfer results from $\mathcal O$
to $\overline{\mathcal O}$.

To describe some of the things that can be obtained in this way,
let $\overline{\Lambda}$ denote a set of representatives
for the equivalence classes of $E \in \Lambda$ viewed up to degree shifts, 
so that $\overline{\Lambda}$ is a complete set of pairwise non-isomorphic 
irreducible admissible $\mathfrak h$-supermodules.
Also let $\widehat{E}$ denote the projective cover of $E \in 
\overline{\Lambda}$ in the category of admissible $\mathfrak h$-supermodules.
We have the objects
$L(E), \De(E)$ and $\nabla(E) \in \overline{\mathcal O}$
obtained from the ones defined before by forgetting the grading.
Intrinsically,
$\De(E) = U(\mathfrak g) \otimes_{U(\mathfrak b)} \widehat{E}$,
$\nabla(E)$ is the largest submodule of 
$\hom_{U(\mathfrak g_{\leq 0})}(U(\mathfrak g), E)$ that belongs to 
$\overline{\mathcal O}$, and
$L(E) = \cosoc_{\mathfrak g} \De(E) \simeq \soc_{\mathfrak g} \nabla(E)$.
In particular, $\{L(E)\}_{E \in \overline{\Lambda}}$ is a complete set
of pairwise non-isomorphic irreducible objects in $\overline{\mathcal O}$.

The notion of a $\De$-flag of an object of $\overline{\mathcal O}$
is defined as before. The multiplicity $(M:\De(E))$
of $\De(E)$ as a subquotient
of a $\De$-flag of an object $M \in \overline{\mathcal O}$ is independent
of the choice of flag. We have that
\begin{equation}
(M:\De(E)) = \dim \hom_{\overline{\mathcal O}}(M, \nabla(E)) / d_E.
\end{equation}
We also note from Corollary~\ref{fp} that summands of objects with finite 
$\De$-flags have finite $\De$-flags.
We can always choose a partial ordering $\preceq$ on 
$\overline{\Lambda}$ such that
\begin{equation*}
[\Delta(F):L(E)] \neq 0 \hbox{ or }
[\nabla(F):L(E)] \neq 0 \Rightarrow E \preceq F.
\end{equation*}
Let $\sim$ be the equivalence relation on $\overline{\Lambda}$
generated by the partial order $\preceq$.
For $\overline{\theta} \in \overline{\mathcal O} / \sim$, 
let
$\overline{\mathcal O}_{\overline{\theta}}$ 
be the full subcategory of $\overline{\mathcal O}$
consisting of the objects $M \in \overline{\mathcal O}$ all of whose
irreducible subquotients are of the form
$L(E)$ for $E \in \overline{\theta}$.
Then Theorem~\ref{blockdec} gives us the block decomposition
of $\overline{\mathcal O}$: 

\begin{Theorem} The functor 
$$
\prod_{\overline{\theta} \in \overline{\Lambda} / \sim}
\overline{\mathcal O}_{\overline{\theta}} \rightarrow \overline{\mathcal O},
\quad
(M_{\overline{\theta}})_{\overline{\theta}} \mapsto 
\bigoplus_{\overline{\theta} \in \overline{\Lambda} / \sim} M_{\overline{\theta}}
$$
is an equivalence of categories.
\end{Theorem}

Next, we use Theorem~\ref{tiltthm} to define
the {\em indecomposable tilting module}
$T(E) \in \overline{\mathcal O}$
for each $E \in \overline{\Lambda}$:

\begin{Theorem}
For each $E \in \overline{\Lambda}$ there exists a unique up to isomorphism
indecomposable object $T(E) \in \overline{\mathcal O}$ such that
\begin{itemize}
\item[(i)] $\ext^1_{\overline{\mathcal O}}(\Delta(F), T(E)) = 0$ for all 
$F \in \overline{\Lambda}$;
\item[(ii)] $T(E)$ admits a $\Delta$-flag starting with $\Delta(E)$ at the 
bottom.
\end{itemize}
\end{Theorem}

Let $\gamma$ be an admissible semi-infinite character for $\mathfrak g$
and construct the semi-regular bimodule $S_\gamma$ as in Lemma~\ref{icky}.
Let $\overline{\mathcal O}^{\Delta}$ be the full subcategory
of $\overline{\mathcal O}$ consisting of the objects that admit a finite
$\Delta$-flag. Then Corollaries~\ref{maineq2} and \ref{dde} give us:

\begin{Theorem}\label{maineq}
The functor $M \mapsto (S_\gamma \otimes_{U(\mathfrak g)} M)^\star$
defines a contravariant 
equivalence of categories $\overline{\mathcal O}^\Delta \rightarrow
\overline{\mathcal O}^\Delta$
under which short exact sequences correspond to short exact sequences
and $\Delta(E)$ maps to
$\Delta(\C_{-\gamma} \otimes E^\#)$ 
for every $E \in \overline{\Lambda}$.
Moreover,
\begin{equation}\label{ds}
(T(E):\Delta(F)) = [\nabla(\C_{-\gamma} \otimes F^\#):L(\C_{-\gamma}\otimes E^\#)]
\end{equation}
for all $E, F \in \overline{\Lambda}$.
\end{Theorem}

Still assuming that (A6) holds, we now
impose some finiteness conditions. First, assume
\begin{itemize}
\item[(A7)] for each $E \in \overline{\Lambda}$,
$\nabla(E)$ has a composition series.
\end{itemize}
Given (A7), it is not hard to show that every object
$M$ in the category $\overline{\mathcal O}^{\text{fin}}$ of 
all {\em finitely generated} admissible $\mathfrak g$-supermodules
that are locally finite dimensional over $\mathfrak b$ has a composition series.
We remark that (A7) holds automatically if the
partial ordering $\preceq$ chosen above has the property that
for each $E \in \overline{\Lambda}$, there are only finitely
many $F \in \overline{\Lambda}$ with $F \preceq E$.
Next assume
\begin{itemize}
\item[(A8)] the category $\overline{\mathcal O}^{\text{fin}}$ has
enough projectives.
\end{itemize}
By Theorem~\ref{pct}(iii), (A8) holds automatically if the
partial ordering $\preceq$ has the property that
for each $E \in \overline{\Lambda}$, there are only finitely
many $F \in \overline{\Lambda}$ with $E \preceq F$.

Using (A7), (A8) and Fitting's lemma, one deduces that each
$L(E)$ has a projective cover denoted $P(E)$ in the category
$\overline{\mathcal O}^{\text{fin}}$. 
Moreover, Theorem \ref{pct}
and Corollary~\ref{bggrec} imply in the present setting that
$P(E)$ has a finite $\Delta$-flag
satisfying BGG reciprocity 
\begin{equation}\label{bgg}
(P(E):\Delta(F)) = [\nabla(F):L(E)]
\end{equation}
for all $E, F \in \overline\Lambda$.
Under the equivalence of categories from Theorem~\ref{maineq},
$P(E)$ gets mapped to $T(\C_{-\gamma} \otimes E^\#)$,
so the tilting modules $T(E)$ also all have {\em finite} $\Delta$-flags,
i.e. they belong to the category $\overline{\mathcal O}^{\text{fin}}$ too.

\section{Examples}

We now give some examples, beginning with the classical ones
to set the scene.

\begin{Example}\rm\label{eg1}
Let $\mathfrak g$ be a finite dimensional semisimple Lie algebra.
Let $\mathfrak t \subset \mathfrak g$ be a maximal toral subalgebra,
and $\Delta\subset \mathfrak t^*$ be a choice of simple roots.
Let $\rho \in \mathfrak t^*$ be half the sum of the corresponding
positive roots. We take the $\Z$-grading on $\mathfrak g$
defined so that $\mathfrak g_{\alpha}$ is in degree $1$
and $\mathfrak g_{-\alpha}$ is in degree $-1$
for each $\alpha \in \Delta$. Clearly this grading is induced
by the adjoint action of some $D \in \mathfrak t$, and
$\mathfrak h := \mathfrak g_{0} = \mathfrak t$.
Taking the group $X$ of admissible weights to be all of $\mathfrak t^*$,
the category $\overline{\mathcal O}^{\text{fin}}$
is exactly the category introduced in \cite{BGG}. 

It is easy to see that our assumptions (A1)--(A6) are all 
satisfied. Moreover, by Harish-Chandra's theorem on central characters,
we can choose the equivalence relation $\sim$ so that the equivalence
classes are the orbits of the finite Weyl group $W$ under the dot action.
Hence the equivalence classes are finite, so 
(A7) and (A8) automatically hold too.
We also note that the usual Verma modules $M(\la)$ for 
$\la \in \mathfrak h^*$ are
the standard modules here,
and their duals under
the duality of \cite[$\S$4, Remark]{BGG} are
the costandard modules.
The indecomposable tilting modules 
$T(\la)$ 
are the modules defined originally by Collingwood and Irving in \cite{CI}.

This setup is generalized to an arbitrary symmetrizable Kac-Moody
algebra in \cite[$\S$7]{so2}, see also \cite{DGK, RC}.
In general, (A7) and (A8) do not hold, so it becomes important to
work in category $\overline{\mathcal O}$ rather than 
$\overline{\mathcal O}^{\text{fin}}$. 
Soergel also discusses certain parabolic analogues.
\end{Example}

\begin{Example}\rm\label{eg2}
In the next two examples, we take $\mathfrak g$ to be the Lie
superalgebra $\mathfrak{gl}(m|n)$.
We recall that $\mathfrak g$ consists of
$(m+n)\times(m+n)$ matrices over $\C$, where we label
rows and columns of such matrices by the ordered index set 
$\{-m,\dots,-1,1,\dots,n\}$. 
Writing $\bar i = \0$ if $i > 0$ and $\1$ if $i < 0$, the parity
of the $ij$-matrix unit $e_{i,j} \in \mathfrak g$ 
is $\bar i + \bar j$, and the superbracket satisfies
$[e_{i,j}, e_{k,l}] = \delta_{j, k} e_{i, l} - (-1)^{(\bar i + \bar j)(\bar k + \bar l)} \delta_{i,l} e_{k,j}.$
The subalgebra $\mathfrak g_{\0}$ of $\mathfrak g$
is isomorphic to $\mathfrak{gl}(m) 
\oplus \mathfrak{gl}(n)$.
We will always take the maximal toral subalgebra $\mathfrak t$
to be the subalgebra consisting of all diagonal matrices, and the group
$X$ of admissible weights to be all of $\mathfrak t^*$.
Let $\delta_{-m},\dots,\delta_{-1},\delta_1,\dots,\delta_n$ be the basis
for $\mathfrak t^*$ dual to the basis $e_{-m,-m},\dots,e_{-1,-1},
e_{1,1},\dots,e_{n,n}$ of 
$\mathfrak t$.

Now there are two natural $\Z$-gradings to consider.
First, we discuss the {\em principal grading} induced
by the adjoint action of the matrix
$D = \operatorname{diag}(m+n,m+n-1,\dots,2,1) \in \mathfrak h$, so
the degree of $e_{i,j}$ is defined by the equation
$[D, e_{i,j}] = \deg(e_{i,j}) e_{i,j}$.
For this grading, $\mathfrak h := \mathfrak g_0$ coincides with 
the subalgebra $\mathfrak t$ of diagonal matrices
and $\mathfrak b := \mathfrak g_{\geq 0}$ is the subalgebra of all
upper triangular matrices. 
Let $\overline{\mathcal O}^{\text{fin}}$ be the resulting
category as in section \ref{var}.
We should check that the assumptions (A1)--(A8) all hold,
the only difficult ones being (A7) and (A8):

\begin{Lemma}\label{eop} 
Every object $M \in \overline{\mathcal O}^{\text{\rm fin}}$ has a composition series, and $\overline{\mathcal O}^{\text{\rm fin}}$ has enough projectives.
\end{Lemma}

\begin{proof}
Let $\mathcal E$ be the category
of all finitely generated $\mathfrak g_{\0}$-supermodules
that are locally finite dimensional over $\mathfrak b_{\0}$ and 
semisimple over $\mathfrak h$.
By the PBW theorem, $U(\mathfrak g)$ is free
of finite rank as a (left or right) $U(\mathfrak g_{\0})$-module.
Hence, the functor
$U(\mathfrak g) \otimes_{U(\mathfrak g_{\0})} ?$ maps objects in
$\mathcal E$ to objects in $\overline{\mathcal O}^{\text{fin}}$, and it is
left adjoint to the natural restriction functor from 
$\overline{\mathcal O}^{\text{fin}}$ to $\mathcal E$.
So it sends projectives to 
projectives. By Example~\ref{eg1}, $\mathcal E$ has
enough projectives, so we deduce that $\overline{\mathcal O}^{\text{fin}}$ 
does too.
Finally, to see that every object $M \in \overline{\mathcal O}^{\text{fin}}$ 
has a composition series,
note that $U(\mathfrak g)$ is Noetherian, so $M$ has a descending filtration
$M = M_0 \geq M_1 \geq \dots$ such that each $M_i / M_{i+1}$ is irreducible.
We just need to show that this filtration stabilizes after finitely many 
terms. But 
every object in $\mathcal E$ has a composition series by Example~\ref{eg1}
so this follows immediately on restricting $M$ to $\mathfrak g_{\0}$.
\end{proof}

The standard modules $\De(\la)$ in this case are the {\em Verma modules}
$M(\lambda) := U(\mathfrak g) \otimes_{U(\mathfrak b)} \C_\lambda$,
where $\C_\lambda$ is the one dimensional $\mathfrak b$-module
with character $\lambda \in \mathfrak h^*$.
The costandard modules $\nabla(\la)$ are the {\em dual Verma modules}
$M(\lambda)^\tau$, where $\tau$ is the duality defined
using the ``supertranspose''
antiautomorphism $e_{i,j} \mapsto (-1)^{\bar i(\bar i + \bar j)} e_{j,i}$
of $\mathfrak g$.
Finally the indecomposable tilting modules are denoted $T(\la)$ and the irreducible
modules are denoted $L(\la)$, for $\la \in \mathfrak h^*$.
Like in Example~\ref{eg1}, an admissible semi-infinite 
character for $\mathfrak g$ with respect to the principal grading
is given by the character
$2\rho$,
where
$\rho = m \delta_{-m}+\dots + 2 \delta_{-2} + \delta_{-1} - \delta_1 - 2\delta_2 - \dots - n \delta_n.$
Now we get from (\ref{ds}) that
\begin{equation}\label{mdt}
(T(\lambda): M(\mu)) = [M(-\mu - 2\rho):L(-\lambda-2\rho)],
\end{equation}
for $\lambda, \mu \in \mathfrak h^*$.
A precise conjecture for these multiplicities in the case
that $\la, \mu$ are integral linear combinations of
the $\delta_i$ can be found in \cite{B1}.

It is interesting to note in this example that both (A7) and (A8)
hold, despite the fact (as seen in \cite{B1}) that 
the partial ordering $\preceq$ of section \ref{var} always
has infinite chains.
\end{Example}

\begin{Example}\rm\label{eg3}
Continuing with $\mathfrak g = \mathfrak{gl}(m|n)$, we now discuss
the second natural $\Z$-grading, namely, the {\em compatible grading}.
This is induced by the adjoint action of the matrix
$D = \operatorname{diag}(1/2,1/2,\dots,1/2;-1/2,-1/2,\dots,-1/2)$.
Note this time that $\mathfrak h := \mathfrak g_0 = \mathfrak g_{\0}$,
and $\mathfrak g_{-1} \oplus \mathfrak g_1 = \mathfrak g_{\1}$.
This time, as is easy to show, the category
$\overline{\mathcal O}^{\text{fin}}$ is precisely the category of
all  finite dimensional
$\mathfrak g$-supermodules that are semisimple over $\mathfrak t$.
The hypothesis (A1)--(A8) are all satisfied, arguing as in Lemma~\ref{eop}
for (A7) and (A8).

Recalling $\mathfrak h \cong \mathfrak{gl}(m)\oplus\mathfrak{gl}(n)$,
the irreducible finite dimensional 
$\mathfrak h$-supermodules 
are parametrized by the set $X^+$ of
{dominant weights}, namely, the
$\la = \la_{-m} \delta_{-m}+\dots+\la_{-1}\delta_{-1}+\la_1\delta_1+\dots+
\la_n \delta_n \in \mathfrak h^*$ with each
$\la_{-m}-\la_{1-m},\dots,\la_{-2}-\la_{-1}, \la_1 - \la_2,\dots,\la_{n-1}-\la_n$ being non-negative integers.
Given $\la \in X^+$, we denote the corresponding standard module
$\Delta(\la)$ instead by $K(\la)$ and call it the {\em Kac module}
of highest weight $\la$, since it
was first defined by Kac in \cite{Kac2}.
The costandard modules are the {\em dual Kac modules} 
$K(\la)^\tau$.
We also write $L(\la)$ for the unique irreducible quotient of $K(\la)$,
$P(\la)$ for its projective cover, and $U(\la)$ for the indecomposable
tilting module
of highest weight $\la$ in this finite dimensional setting.
By (\ref{bgg}), $P(\la)$ has a finite
Kac flag with $K(\la)$ at the top, satisfying the BGG reciprocity
\begin{equation}
(P(\la):K(\mu)) = [K(\mu):L(\la)],
\end{equation}
as was also proved in \cite[Proposition 2.5]{Zou}.

Now let
$\beta = n(\delta_{-m}+\dots+\delta_{-1}) - m (\delta_1+\dots+\delta_n)$
be the sum of the positive odd roots.
It is easy to check that the unique $1$-dimensional representation
$\gamma:\mathfrak h\rightarrow\C$ of
weight $-\beta$ is an admissible semi-infinite character for $\mathfrak g$ with
respect to the compatible grading.
In fact in this case, there is an
even isomorphism of $U(\mathfrak g), U(\mathfrak g)$-bimodules
between the semi-regular bimodule $S_{\gamma}$
from Lemma~\ref{icky} and the regular bimodule $\Pi^{mn} U(\mathfrak g)$.
In the notation of Lemma~\ref{icky}, 
an isomorphism maps $1 \in U(\mathfrak g)$ to the
element $\iota(\delta) \in S_\gamma$, where $\delta \in U(\mathfrak n)^{\circledast}$ is the function mapping $\prod_{-m \leq i < 0 < j \leq n} e_{j,i}$
to $1$
(product taken in some fixed order)
and all other monomials in the $e_{j,i}$ of strictly smaller length to $0$.
So in this case the duality in Theorem~\ref{maineq} is (up to 
parity change and degree shift) 
just the usual duality $*$ on finite dimensional
$\mathfrak g$-supermodules. In particular, 
\begin{align}
K(\beta - w_0\la)^* \cong K(\la),\label{kdual}\\
P(\beta-w_0\la)^* \cong U(\la),\label{pdual}
\end{align}
where $w_0$ denotes the longest element of the Weyl group
$W \cong S_m \times S_n$ of $\mathfrak h$ 
acting on $\mathfrak t^*$ in the obvious way.
The statement (\ref{ds}) says
\begin{equation}\label{kdt} 
(U(\lambda): K(\mu)) = [K(\beta-w_0 \mu):
L(\beta-w_0 \lambda)],
\end{equation}
for $\lambda, \mu \in X^+$.
The numbers on the left hand side of this equation are
computed in \cite{B1}.
\end{Example}

\begin{Example}\rm\label{eg4}
In the final example, we take $\mathfrak g = \mathfrak{q}(n)$.
Thus, $\mathfrak g$ is the subalgebra of $\mathfrak{gl}(n|n)$
consisting of all matrices of the form
$\left(
\begin{array}{l|l}X&Y\\\hline Y&X\end{array}
\right)$.
For $1 \leq i,j\leq n$, 
we will let $e_{i,j}$ resp. $e_{i,j}'$ denote the even resp.
odd matrix unit, i.e. the matrix of the above form
with the $ij$-entry of $X$ resp. $Y$ equal to $1$ and all other entries
equal to zero. The $\Z$-grading on $\mathfrak g$ is defined by
$\deg(e_{i,j}) = \deg(e_{i,j}') = (j - i)$.
For this grading, $\mathfrak h := \mathfrak g_0$ 
is spanned by $\{e_{i,i}, e_{i,i}'\:|\:1 \leq i \leq n\}$,
and $\mathfrak b := \mathfrak g_{\geq 0}$ is 
spanned by $\{e_{i,j}, e_{i,j}'\:|\:1 \leq i \leq j \leq n\}$.
We also let $\mathfrak t = \mathfrak h_{\0}$ and take the group $X$
of admissible weights to be all of $\mathfrak t^*$.

As explained in \cite[$\S$3]{penkov},
the finite dimensional irreducible $\mathfrak h$-supermodules
are parametrized by the set $\mathfrak t^*$.
For $\la \in \mathfrak t^*$, we write $\u(\la)$ for the corresponding
irreducible $\mathfrak h$-supermodule. It is constructed in
\cite{penkov} as a certain Clifford module,
of dimension a power of $2$. The assumption (A5)
can be checked from this construction and the fact that Clifford algebras
are symmetric: one gets that
\begin{equation}
\widehat{\u(\la)}^* \cong \widehat{\u(-\la)}.
\end{equation}
The remaining assumptions (A1)--(A4) and (A6) are easy, 
and one argues like in Lemma~\ref{eop} to verify (A7) and (A8).

So now we can consider the category $\overline{\mathcal O}^{\text{fin}}$
as in section~\ref{var}.
Let
\begin{align*}
M(\la) := U(\mathfrak g) \otimes_{U(\mathfrak b)} \u(\lambda),\qquad
N(\la) := U(\mathfrak g) \otimes_{U(\mathfrak b)} \widehat{\u(\la)},
\end{align*}
for each $\la \in \mathfrak t^*$.
Then $N(\la)$ is the standard module $\De(\la)$
in $\overline{\mathcal O}^{\text{fin}}$,
while $M(\la)$
is dual to the costandard module $\nabla(\la)$  
under the duality $\tau$ induced by the (unsigned)
antiautomorphism 
$$
\left(
\begin{array}{l|l}X&Y\\\hline Y&X\end{array}
\right)
\mapsto\left(
\begin{array}{l|l}X^T&Y^T\\\hline Y^T&X^T\end{array}
\right).
$$
One checks that the trivial character $0:\mathfrak h \rightarrow \C$
is an admissible semi-infinite character for $\mathfrak g$.
So, 
writing $T(\la)$ resp. $L(\la)$ for the indecomposable
tilting module resp. the irreducible
module
corresponding to $\la \in \mathfrak t^*$, (\ref{ds}) shows that
\begin{equation}\label{mc1}
(T(\la):N(\mu)) = [M(-\mu):L(-\la)].
\end{equation}
A precise conjecture for these decomposition numbers in case
$\la,\mu$ are integral weights is
formulated in \cite{B2}.
\end{Example}

\end{document}